# SHORT NOTE ON ADDITIVE SEQUENCES AND ON RECURSIVE PROCESSES

Andrei Vieru

**Abstract**

Simple methods permit to generalize the concepts of iteration and of recursive processes. We shall see briefly on several examples what these methods generate. In additive sequences, we shall encounter not only the golden or the silver ratio, but a *dense set* of ratio limits that corresponds to an infinity of conceivable recursive additive rules. We shall show that some of these limits have nice properties. Identities involving Fibonacci and Lucas sequences will be viewed as special cases of more general identities. We shall show that some properties of the Pascal Triangle belong also to other similar objects. In Dynamical Systems and Chaos Theory we shall encounter *weird orbits*, whose order is higher than the number of its distinct elements and, beyond the *chaos point*, a rather unexpected belated convergence to 0, after a pseudo chaotic behaviour during as many terms as one may wish. Time and again, we shall find here the Feigenbaum constant. In Formal Grammars we shall see that recursive rules applied to concatenation are sometimes equivalent to formal grammars although generally more restrictive.

## 1. RECURRENCE RELATIONS IN CHAOS THEORY

### 1.1. INTERESTING CASES OF GENERALIZED ITERATION

We already have studied what we call generalized iteration[1]; we showed that generalized iteration of the first and of the second kind don't imply – in spite of common ideas – the same number of dimensions of the domain and of the range.

---

[1] see our article '*Generalized iteration, catastrophes, generalized Sharkovsky's Ordering*' arXiv:0801.3755v2 **[math.DS]**

For example, we already studied the recursive model

$$u_n = F_a(u_{n-2}, u_{n-1}) \tag{1°}$$

bringing to the fore some interesting results, including the appearance of chaos, of the Feigenbaum constant (when parameter $a$ increases) and of what we call 'generalized Sharkovsky's ordering', not to speak about some discontinuities gotten with continuous ingredients.

Consider the following sequence, where at each step all occurrences of $x$ are replaced by $y$ and all occurrences of $y$ are replaced by $F_a(x, y)$:

$x$,

$y$,

$F_a(x, y)$,

$F_a(y, F_a(x, y))$,

$F_a(F_a(x, y), F_a(y, F_a(x, y)))$,

$F_a(F_a(y, F_a(x, y)), F_a(F_a(x, y), F_a(y, F_a(x, y))))$

etc. where $F_a(x, y) = ax(1-x)y(1-y)$

As one can see, the sequence is based on the recursive rule

$$\forall n > 2 \quad u_n = F(u_{n-2}, u_{n-1}) \tag{2°}$$

This may be also written as a set of two rules of replacement at each step of the recursive process of all occurrences of $x$ and $y$:

$$x \rightarrow y, \ y \rightarrow F_a(x, y) \tag{2bis°}$$

F-iterating a function $t = f(x, y)$ we may take, arbitrarily, more than 2 initial values.

One can try to see what happens if we fix a set of recursive rules of the following type (to be applied at each step to all occurrences of the involved variables):

$$x \rightarrow y, y \rightarrow z, z \rightarrow F(x, y) \tag{3°}$$

Starting from $x$ they yield the following sequence:

$x$,

$y$,

$z$,

$F(x, y)$

$F(y, z)$

$F(z, F(x, y))$

$F(F(x, y), F(y, z))$

$F(F(y, z), F(z, F(x, y)))$

etc

This sequence is based on the recursive rule $\forall n>3\ u_n = F(u_{n-3}, u_{n-2})$     **(4°)**

Another interesting set of rules is $x \to y, y \to z, z \to F(x, z)$     **(5°)**

It yields the sequence:

$x,$

$y,$

$z,$

$F(x, z)$

$F(y, F(x, z))$

$F(z, F(y, F(x, z)))$

$F(F(x, z), F(z, F(y, F(x, z))))$

etc.

which is based on the recursive rule

$\forall n>3\ u_n = F(u_{n-3}, u_{n-1})$     **(6°)**

It easy to see that if, for instance, we chose $F_a(x, y) = ax(1-x)y(1-y)$, both produce chaos but not exactly through an usual bifurcation process.

For example, the recursive rule **(6°)** applied to $F_a(x, y)=ax(1-x)y(1-y)$ produce first a normal doubling-period (at $a \approx 10.415$, we have a bifurcation and then a period-two orbit), then at about $a \approx 12.53$, we switch directly to a *cycle* of order 8 in which the first point equals the 5th point, the 2nd point equals the 8th point and the 4th equals the 6th. The cycle has the form $(a, b, c, d, a, d, e, b)$.

Sets of rules such as **(2°bis), (3°)** and **(5°)** might be called ***iterating operators.***

*Stable* cycles with periodicity higher than the number of their own *distinct* elements are, if we make no mistake, unknown in other contexts.

This kind of cycles, along with the context they appear in, leads to the following definition, which generalizes the definition of periodic orbits.

## 1.2. PERIODIC CAROUSELS

Let $\Omega$ be an iterating operator that transforms a map $g_i : \mathbf{T}^n \to \mathbf{T}^m$ ($m \le n$) into a map $\Omega(g_i): \mathbf{T}^n \to \mathbf{T}^m$, ($m \ge n$). (We'll write $\Omega^0(g_i) = g_i$) For a map $g_0$ a **periodic carousel** with **least** period $k$ is the finite sequence of $k$ (not necessarily distinct) points in $\mathbf{T}^m$ $\{p_j\}_{0 \le j \le k-1}$ with $p_j = (\Omega^{j-1}(g_0))(p_0)$ ($0 < j < k$), with $(\Omega^{k-1}(g_0))(p_0) = p_0$, with, for any non-negative integers $l$ and $p$,

$$l \equiv p \ (\text{mod } k) \Rightarrow (\Omega^l(g_0))(p_0) = (\Omega^p(g_0))(p_0)$$

and such that for any $k' < k$ there is at least a pair of distinct integers $l$ and $p$ with $l \equiv p \ (\text{mod } k')$ and $(\Omega^l(g_0))(p_0) \ne (\Omega^p(g_0))(p_0)$

Any periodic orbit for any map $g$ is also a periodic carousel if we let $\Omega^n(g) = g$ for any $g$ and any $n$ (or if $\Omega$ is idempotent).

Under generalized iteration the concept of *periodic carousel* applies to maps $g_i : \mathbf{R}^n \to \mathbf{R}^m$, $n \ge m$

A periodic carousel with period $k$ is said to be *normal* if the $k$ points that form it are all distinct. (All periodic orbits are normal: the characteristics of iteration of single-valued functions guarantee their normality.)

## 1.3. WEIRD PERIODIC CAROUSELS

A periodic carousel with period $k$ is said to be *weird* if at least two of its points are identical. In other words a weird carousel crosses itself.

**CONJECTURE**

Weird periodic carousels occur if and only if the iterated mapping is symmetric with respect to at least one of the axis of its domain.

## 1.4. HOW GENERAL THE FEIGENBAUM CONSTANT IS?

Then, at $a \approx 13.27$, we'll have a bifurcation toward a normal orbit of order 16 with 16 distinct elements. The next bifurcations points are: (about) 13.417…, 13.4515…, 13.4593…, 13.46102…, 13.46139… So, after some computations, we see that it is likely the Feigenbaum constant appears here too.

Amazingly, when *a* approaches values beyond 15.595, the whole sequence converges to 0. Moreover, the converging to 0 sequences may 'seem chaotic' from their first term until some term **whose rank may be as high as one wishes**. For *a*>15.32, the recursive rule **(4°)** also produces convergence to 0, in an even stranger way: **it is not at all clear in which way the 'speed' of this convergence (i.e. in which way the length of the initial pseudo-chaotic subsequence) depends on the parameter *a* value.**

After periods of chaos we 'sink' again into order, i.e. in a stable periodic orbit of order 7 (when $a \in [14.6, 14.9]$)

(In the experiment described above, the initial values where chosen as follows: $x = 0.6$; $y = 0.7$; $z = 0.8$.)

The recursive rule **(4°)** engenders a somewhat different behavior. We have first a split of the fixed point into an order 5 periodic orbit[2] (for *a* slightly higher than 11). Then, chaos is reached extremely quickly through doubling period, when *a* grows from 13.165… to 13.17… (One can try to determine if here the *Feigenbaum constant* appears again or not.)

Once again, amazingly, when *a* > 15.3… the sequences converge to 0. Moreover, the converging to 0 sequences may seem chaotic from their first term until some term whose rank may be as high as one wishes.

## 2. AN EXCURSION INTO NUMBER THEORY

Everyone knows what Fibonacci numbers are. Someone has had the idea to generalize the concept of additive sequence. The so-called sequences of *k*-bonacci numbers – whose every term of rank > *k* is the sum of the *k* previous terms – were studied. There is always a limit the ratio of two consecutive terms is converging to. This limit is always a root of the equation $x^n - x^{n-1} - x^{n-2} - \ldots - x - 1 = 0$, namely the root situated in the interval ]1, 2[ (such a root always exists.) The greater the *n* value, nearer the root will be situated to the right extremity of the open interval ]1, 2[.

The sequence of these limits converges to 2. There is a beautiful formula:

---

[2] One can note that the *weird orbit* with period 8 obtained with the recursive rules **(5°) - 6°)** (see above) has exactly 5 distinct elements.

$$\frac{\ln\left[(2-x)^{-1}\right]}{\ln(x)} = 2$$

that establishes a relation between such a limit – $x$ – and the number $n$ of preceding terms whose sum determines each of the following ones.

## 2.1. RECURRENCE EQUATIONS IN NUMBER THEORY

Additive sequences can be conceived in a broader way. Any additive sequence is, in fact, a special case of what we call 'generalized iteration[3]' of a function $f$ : $R^n \to R$ (in the Fibonacci numbers case, we have $n=2$ and $f(p, m)=p+m$; in the Tribonacci numbers case, we have $n=3$ and $f(l, m, p)=l+m+p$, etc.)

The Padovan numbers constitute a well-known example which shows that the number of the arbitrarily chosen first terms of an additive sequence don't have to equal the number of terms whose sum defines the general term of the considered additive sequence.

It is based on the following recurrence equation:

$\forall n>3$  $u_n = u_{n-2} + u_{n-3}$ **(I°)**

It is easy to see that

$\lim_{n \to \infty} u_{n+1}/u_n = 1.3247179572447\ldots$

which is one of the roots of the equation $x^3-x-1=0$, namely the real root, situated in the interval ]1, 2[.

The number 1.32471… , which will be designated hereafter by $\phi_4$ – because of the logarithmic equation[4] **(II°)** – is also a root of the quintic equation[5] $x^5-x^2-x-1=0$, of the quintic equation[6] **$x^5-x^4-1=0$**, of the sixtic equation $x^6-x^4-x-1=0$, of the equation $x^7-x^4-x^2-x-1=0$, of the equation $x^8-x^7-x-1=0$, of the equation $x^8-x^4-x^3-x^2-x-1=0$, of the equation $x^9-2x^4-x^3-x^2-x-1=0$, of the equation $x^{10}-x^9-x^4-1=0$ of the equation $x^{10}-x^7-x^6-x^4-1=0$, of the equation $x^{10}-x^5-2x^4-x^3-x^2-x-1=0$ etc. and of the equation

$$\frac{\ln\left[(x-1)^{-1}\right]}{\ln(x)} = 4 \qquad \text{(II°)}$$

---

[3] See our article '*Generalized iteration, catastrophes, generalized Sharkovsky's Ordering*' arXiv:0801.3755v2 **[math.DS]**

[4] equations of the form $x^{n+1} - x^n - 1 = 0$ are equivalent to $\text{Log}[1/(x-1)]/\text{Log}(x)=n$. We'll keep mentioning the logarithmic form, because, at least for the Pisot numbers when $n=2, 3$ or $4$ there are some nice equalities.

[5] It is therefore also the ratio limit of the additive sequence based on the recursive rule $\forall n>5$  $u_n = u_{n-3} + u_{n-4} + u_{n-5}$ and of an infinity of other recursive rules that may easily be deduced from the quoted equations (see below).

[6] the equations written in **bold** belong to a recursive series of equations that have a beautiful property (see below, Statement 1).

Replacing the recurrence relation **(I°)** by

$$\forall n>3 \ u_n=u_{n-1}+u_{n-3} \qquad \text{(III°)}$$

we'll find $\lim_{n\to\infty} u_{n+1}/u_n = 1.465571231876768\ldots$, which is the root of the equation $x^3-x^2-1=0$ situated in the interval $]1, 2[$.

The number $\phi_2=1.46557123\ldots$ is also a root of the equation[7] $x^4-x^2-x-1=0$, of the equation $x^5-2x^2-x-1=0$, etc., and of the equation

$$\frac{\ln\left[(x-1)^{-1}\right]}{\ln(x)} = 2$$

$\phi_3=1.3802775\ldots$ is a root of the of the equation $x^4-x^3-1=0$ (and therefore the limit ratio of every additive sequence whose rule is $\forall n>4 \ u_n=u_{n-1}+u_{n-4}$), of the quintic equation $x^5-x^3-x-1=0$, of the sixtic equation $x^6-x^3-x^2-x-1=0$, of the equation $x^9-x^7-x^5-x^3-x-1=0$ and of the equation

$$\frac{\ln\left[(x-1)^{-1}\right]}{\ln(x)} = 3$$

Again, as for the other previously examined numbers, the algebraic equations whose root this number is, enables us to formulate the recurrence relations that generates sequences with a ratio limit which equals it.

$\phi_5=1.2851990332\ldots$ is a root of the of the equation $x^6-x^5-1=0$, of the equation $x^7-x^5-x-1=0$, of the equation $x^8-x^5-x^2-x-1=0$, of the equation $x^9-x^5-x^3-x^2-x-1=0$, of the equation $x^{10}-x^5-x^4-x^3-x^2-x-1=0$, etc. and of the equation

$$\frac{\ln\left[(x-1)^{-1}\right]}{\ln(x)} = 5$$

---

[7] It is therefore also the ratio limit of the additive sequence based on the recursive rule $\forall n>4 \ u_n=u_{n-2}+u_{n-3}+u_{n-4}$

$\phi_1 = \phi = 1.618033988\ldots$ is a root of the equation $x^2-x-1=0$, but also of the equations $x^3-2x-1=0$, $x^4-x^2-2x-1=0$, $x^5-x^3-x^2-2x-1=0$, etc.

It is easy to see that every Fibonacci number, as well as every Lucas number[8], may be obtained in infinitely many ways from some previous Fibonacci (respectively, Lucas) numbers. To mention only two of them,

$\forall n>3\ F_n = 2F_{n-2} + F_{n-3}$, $\forall n>4\ F_n = F_{n-2} + 2F_{n-3} + F_{n-4}$ (etc.)

These recursive rules can be almost trivially deduced from the main recursive rule (by mere substitution).

## STATEMENT 1

Let $\phi_k$ be a real in $]0, 1[$ that satisfies the equation

$$\frac{\ln\left[(x-1)^{-1}\right]}{\ln(x)} = k \quad (k \in N) \qquad \textbf{(IV°)}$$

$\forall k \geq 1\ \phi_k$ is root (the only one in $]0, 1[$) of the equation $x^{k+1}-x^k-1=0$, of the equation $x^{k+2}-x^k-x-1=0$, and of any equation of an infinite series $S_k$ of equations, each of whom is obtained from the previous one replacing the highest order term $x^{k+m}$ by the expression $x^{k+m+1} - x^m$

If we write these algebraic equations as

$\Psi_{k,m}(x) = 0 \quad (k, m) \in N \times N$

then, for every given $k$, will have:

$\forall i\ \forall j\quad i < j \Leftrightarrow \Psi'_{k,i}(\phi_k) < \Psi'_{k,j}(\phi_k)$

Moreover

For any $k$ and for any $m > 1$ $\quad \dfrac{\Psi'_{k,m+1}(\phi_k) - \Psi'_{k,m}(\phi_k)}{\Psi'_{k,m}(\phi_k) - \Psi'_{k,m-1}(\phi_k)} = \phi_k \qquad \textbf{(V°)}$

One can ask the question whether the sufficient condition **(IV°)** is necessary to make **(V°)** hold or whether it may be relaxed. One can also ask in which way are formed, for $k > 2$, the equations that have $\phi_k$ as one of their roots and which do not belong to the recursive series of equations described in the statement.

---

[8] And, more generally, of any additive sequence

**STATEMENT 2**

If $k$ is even, then all equations in the series of equations $S_k$ have two and only two real roots, namely $-1$ and $\phi_k$.

If $k$ is odd, then all equations of odd degree in $S_k$ have three and only three real roots, namely $-1$, $\phi_k$ and a number $\sigma_k \in ]-1, 0[$, while all equations in $S_k$ of even degree have only two real roots, namely $\phi_k$ and $\sigma_k$. (In particular all equations of the series $S_1$ have $-1/\phi_1 = -0.61803\ldots$ as one of their roots.)

Moreover, for any odd $k$ and for any integer $m$ greater than 1

$$\frac{\Psi'_{k,2m+4}(\sigma_k) - \Psi'_{k,2m+2}(\sigma_k)}{\Psi'_{k,2m+2}(\sigma_k) - \Psi'_{k,2m}(\sigma_k)} = \sigma_k^2 \qquad \textbf{(VI°)}$$

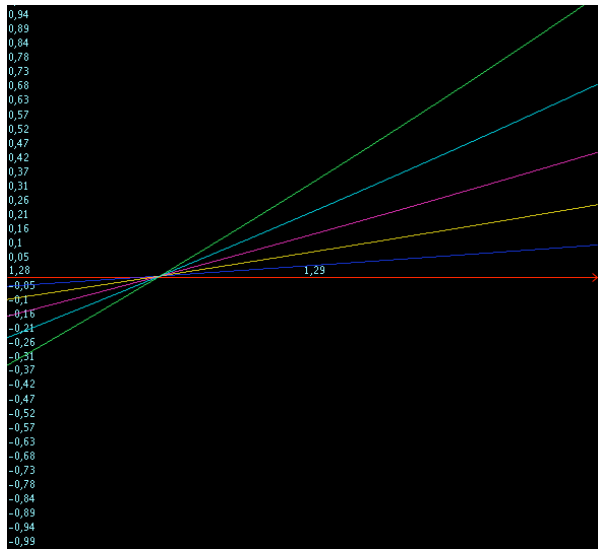

$y=x^6-x^5-1$    $y=x^7-x^5-x-1$    $y=x^8-x^5-x^2-x-1$   (picture 1)
$y=x^9-x^5-x^3-x^2-x-1$    $y=x^{10}-x^5-x^4-x^3-x^2-x-1$

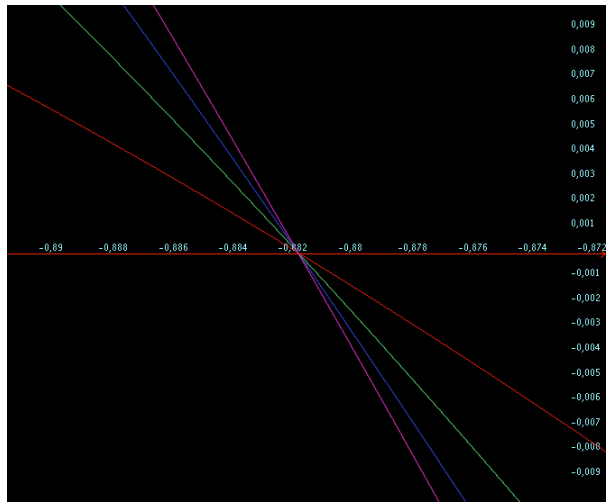

$y=x^7-x^5-x-1$   $y=x^9-x^5-x^3-x^2-x-1$   $y=x^{11}-2x^5-x^4-x^3-x^2-x-1$   (picture 2)
$y=x^{13}-x^7-x^6-2x^5-x^4-x^3-x^2-x-1$

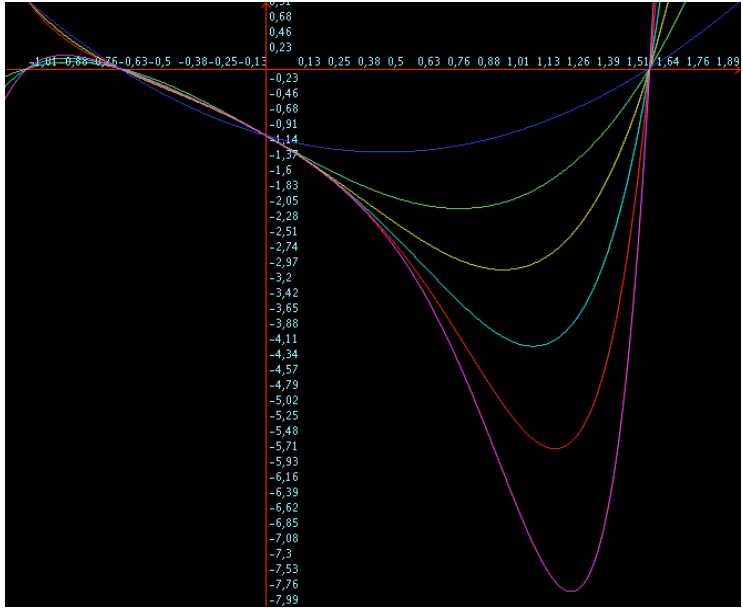

$y=x^2-x-1$  $\qquad y=x^3-2x-1$  $\qquad y=x^4-x^2-2x-1$ (picture 3)

$y=x^5-x^3-x^2-2x-1$  $\qquad y=x^6-x^4-x^3-x^2-2x-1$  $\quad y=x^7-x^5-x^4-x^3-x^2-2x-1$

Remarkably $\phi_1$ – the golden ratio – is a solution of both

$$\frac{\ln\left[(2-x)^{-1}\right]}{\ln(x)} = 2 \quad \text{and} \quad \frac{\ln\left[(x-1)^{-1}\right]}{\ln(x)} = 1$$

$\phi_1^2$ is a solution of the equation $\dfrac{\text{Log}\left[(3-x)^{-1}\right]}{\text{Log}(x)} = 1$

$\phi_1^3$ is a solution of the equation $\dfrac{\text{Log}\left[(x-4)^{-1}\right]}{\text{Log}(x)} = 1$

$\phi_1^4$ is a solution of the equation $\dfrac{\text{Log}\left[(7-x)^{-1}\right]}{\text{Log}(x)} = 1$

$\phi_1^n$ is a solution of the equation $\dfrac{\text{Log}\left\{\left[(-1)^{n+1}x + (-1)^n L_n\right]^{-1}\right\}}{\text{Log}(x)} = 1$, where $L_n$ is the

$n$ - th Lucas number.

It follows that

$$\frac{1}{(-1)^{n+1}\phi_1^n + (-1)^n L_n} = \phi_1^n \qquad \textbf{(VII°)}$$

which is a stronger formula[9] than the identity

$L_n = [\phi^n]$ (given by http://mathworld.wolfram.com/LucasNumber.html, where $[x]$ denotes the *nint* – nearest integer – function)

---

[9] It can be derived from the closed form of the Lucas numbers.

More generally, if $S_k$ is the positive solution of the equation $x^2-kx-1=0$ ($k \geq 1$) – in other words, if $S_k$ is a silver mean – then it can be represented by the continued fraction[10] $[k; k, k, k, \ldots]$ and we'll always have

$$\frac{\ln\left\{\left[(-1)^{n+1} S_k^n + (-1)^n L_k(n)\right]^{-1}\right\}}{\ln(S_k)} = n$$

(where $L_k(n)$ is the $n$-th term of the additive sequence $k$, $k^2+2$, $k+k(k^2+2)$, $k^2+2+k[k+k(k^2+2)],\ldots$etc.)

Generalizing **(VII°)** we get $\quad \dfrac{1}{(-1)^{n+1} S_k^n + (-1)^n L_k(n)} = S_k^n$

## 2.2. GENERALIZED LUCAS-PELL AND FIBONACCI-PELL SEQUENCES

Let us consider the additive sequences defined by $u_1=1$, $u_2=j$ ($j \in \mathbf{N}$) and

$$\forall n>2 \quad u_n = u_{n-2} + ju_{n-1} \qquad \textbf{(II)}$$

It is well known that the limit ratio of two consecutive terms of such sequences equals $\frac{1}{2}\left(j + \sqrt{j^2+4}\right)$, which is a silver mean. Let's call such sequences Fibonacci-Pell sequences of order $j$, or $j$-Fibonacci sequences (made of $j$-Fibonacci Numbers, designated by $F_{(j, n)}$).

Let us consider additive sequences defined by $u_1=j$, $u_2=j^2+2$ and

$$\forall n>2 \quad u_n = u_{n-2} + ju_{n-1} \qquad \textbf{(II)}$$

Let's call these sequences Lucas-Pell sequences of order $j$, or $j$-Lucas sequences (made of $j$-Lucas Numbers, designated by $L_{(j, n)}$).

It turns out that many identities that establish a relationship between Fibonacci and Lucas numbers (of order one) hold for any order.

---

[10] There is a celebrated presentation of the golden ratio in a nested radical form. It can be generalized: the positive root of the equation $x^2-(2k-1)x-1=0$ ($k \geq 1$), equals

$$k-1+\sqrt{k^2-k+1+\sqrt{k^2-k+1+\sqrt{k^2-k+1+\sqrt{k^2-k+1+\sqrt{k^2-k+1}}}}}\ldots$$

Soon we'll see that it can be also generalized in another way.

Thus, the well known identity $F_{2n} = L_n F_n$ is a special case of the identity

$$F_{(j, 2n)} = L_{(j, n)} F_{(j, n)} \tag{III}$$

The well known identity $F_{m+n} = \frac{1}{2}(F_m L_n + L_m F_n)$ is a special case of the identity

$$F_{(j, m+n)} = \frac{1}{2}\left(F_{(j, m)} L_{(j, n)} + L_{(j, m)} F_{(j, n)}\right) \tag{IV}$$

The well known identity $F_m L_n = F_{m+n} + (-1)^n F_{m-n}$ is a special case of the identity

$$F_{(j, m)} L_{(j, n)} = F_{(j, m+n)} + (-1)^n F_{(j, m-n)} \tag{V}$$

The well known identity $F_n = (L_{n-1} + L_{n+1})/5$ is a special case of the identity

$$F_{(j, m+n)} = \frac{L_{(j, n-1)} + L_{(j, n+1)}}{j^2 + 4} \tag{VI}$$

The well known identity $L_n^2 + 5F_n^2 = 4(-1)^n$ is a special case of the identity

$$L_{(j, n)}^2 + (j^2 + 4) F_{(j, n)}^2 = 4 \cdot (-1)^n \tag{VII}$$

A. Mihailov's product expansions
$$F_m F_n = [L_{m+n} - (-1)^n L_{m-n}]/5$$
and $$F_m L_n = F_{m+n} + (-1)^n F_{m-n}$$

are special cases of the product expansions

$$F_{(j, m)} F_{(j, n)} = [L_{(j, m+n)} - (-1)^n L_{(j, m-n)}]/(j^2 + 4) \tag{VIII}$$

and
$$F_{(j, m)} L_{(j, n)} = F_{(j, m+n)} + (-1)^n F_{(j, m-n)} \tag{IX}$$

The square expansion $F_n^2 = [L_{2n} - 2(-1)^n]/5$ is a special case of

$$F_{(j, n)}^2 = \frac{L_{(j, 2n)} - 2 \cdot (-1)^n}{j^2 + 4} \tag{X}$$

$F_{2n} = F_{n+1}^2 - F_{n-1}^2$ should be understood as a special case of

$$F_{(j,\,2n)} = \frac{F_{(j,\,n+1)}^2 - F_{(j,\,n-1)}^2}{j} \qquad \text{(XI)}$$

Other general identities:

$$\begin{aligned} F_{(j,\,2n)} &= F_{(j,\,n)}(F_{(j,\,n+1)} + F_{(j,\,n-1)}) \\ &= F_{(j,\,n)}(jF_{(j,\,n+1)} + 2F_{(j,\,n-1)}) \\ &= F_{(j,\,n)}(2F_{(j,\,n+1)} - jF_{(j,\,n)}) \end{aligned} \qquad \text{(XII)}$$

$$F_{(j,\,3n)} = \frac{F_{(j,\,n+1)}^3 - jF_{(j,\,n)}^3 - F_{(j,\,n-1)}^3}{j} \qquad \text{(XIII)}$$

Catalan's and Cassini's identities[11] hold without changes for every *j*, **while the Gelin-Cesàro identity does not.** The following one holds for every *j, n* and *k*:

$$\sum_{k=1}^{n} F_{(j,\,k)}^2 = j^{-1} F_{(j,\,n)} F_{(j,\,n+1)} \qquad \text{(XIV)}$$

Johnson gives a general identity $F_a F_b - F_c F_d = (-1)^r (F_{a-r} F_{b-r} - F_{c-r} F_{d-r})$ for arbitrary integers and with *a+b=c+d*, from which many other identities follow as special cases.

It also holds for any *j*:

$$F_{(j,\,a)}F_{(j,\,b)} - F_{(j,\,c)}F_{(j,\,d)} = (-1)^r (F_{(j,\,a-r)}F_{(j,\,b-r)} - F_{(j,\,c-r)}F_{(j,\,d-r)}) \qquad \text{(XV)}$$

and, as a matter of fact, may be completed as follows:

If min{*a, b, c, d*} is odd and min{*a, b, c, d*} ≠ *a* ≠ max{*a, b, c, d*} or
if min{*a, b, c, d*} is even and min{*a, b, c, d*} ≠ *c* ≠ max{*a, b, c, d*}, then

$$F_{(j,\,a)}F_{(j,\,b)} - F_{(j,\,c)}F_{(j,\,d)} = (-1)^{|c-b|+|c-a|} F_{(j,\,|c-b|)} F_{(j,\,|c-a|)} \qquad \text{(XVI a)}$$

---

[11] $F_n^2 - F_{n-r}F_{n+r} = (-1)^{n-r}F_r^2$, respectively $F_{n-1}F_{n+1} - F_n = (-1)^n$

If min{a, b, c, d} is even and min{a, b, c, d} ≠ a ≠ max{a, b, c, d} or
if min{a, b, c, d} is odd and min{a, b, c, d} ≠ c ≠ max{a, b, c, d}, then

$$F_{(j,\,a)}F_{(j,\,b)} - F_{(j,\,c)}F_{(j,\,d)} = (-1)^{|c-b|+|c-a|+1} F_{(j,\,|c-b|)}F_{(j,\,|c-a|)} \quad \textbf{(XVI b)}$$

Using **(XV)**, **(XVI a)** and **(XVI b)**, one can prove that

$$F_{(j,\,n)}^{4} - F_{(j,\,n+1)}F_{(j,\,n-1)}F_{(j,\,n+2)}F_{(j,\,n-2)} = (-1)^{n}(j^{2}-1)F_{(j,\,n)}^{2} + j^{2} \quad \textbf{(XVII)}$$

Setting $j = 1$, one obtains **the Gelin-Cesàro identity** as a special case.

Many other identities involving Fibonacci and Lucas numbers might be generalized. Of course, many algebraic, combinatorial and number theory problems (that involve primes, powers, triangular numbers, etc.) already studied in the context of 'simple' Fibonacci numbers, deserve to be studied in the larger context of $j$-Fibonacci sequences and $j$-Lucas numbers.

## 2.3. SOME PROPERTIES OF THE $\phi_n$ NUMBERS. DIRECTED GRAPHS AND NESTED RADICAL FORMS

Turning back to the $\phi_n$ numbers, we find the following relations:
$\phi_2^2 = 1.46557123\ldots^2 = 2.14789903\ldots$ is a solution of the equation

$$\frac{\ln\left[(x-2)^{-1}\right]}{\ln(x)} = \frac{5}{2} \quad (\alpha)$$

$\phi_2^3 = 1.46557123\ldots^3 = \phi_2^2 + 1 = 3.14789903\ldots$ is a solution of the equation

$$\frac{\ln\left[(x-3)^{-1}\right]}{\ln(x)} = \frac{5}{3} \quad (\beta)$$

Besides,

$$\frac{\ln\left[(\phi_2^4 - 2)^{-1}\right] - \ln\left[(5 - \phi_2^4)^{-1}\right]}{\ln(\phi_2^4)} = -\frac{5}{4} \quad (\chi)$$

$\phi_3{}^2==1.3802775\ldots{}^2=1.9051661677540189\ldots$ is a solution of the equation

$$\frac{\ln\left[(2-x)^{-1}\right]-\ln\left[(x-1)^{-1}\right]}{\ln(x)}=\frac{7}{2} \qquad (\delta)$$

$\phi_4{}^2=1.32471795\ldots{}^2=1.754877\ldots$ is a solution of the equation $x^8-x^7-x^6-x^3-x^2-1=0$, of the equation $x^4-x^3-x^2-1=0$ and of the equation

$$\frac{\ln\left[(x-1)^{-1}\right]}{\ln(x)}=\frac{1}{2} \qquad (\varepsilon)$$

$\phi_4{}^3=1.32471795\ldots+1=2.3247195\ldots$ is a solution of the equation

$$\frac{\ln\left[(x-2)^{-1}\right]}{\ln(x)}=\frac{4}{3} \qquad (\gamma)$$

$\phi_4{}^4=3.0795956\ldots$ is a solution of the equation

$$\frac{\ln\left[(x-3)^{-1}\right]}{\ln(x)}=\frac{9}{4} \qquad (\rho)$$

$\phi_4{}^5=\phi_4{}^4+1=4.0795956\ldots$ is a solution of the equation

$$\frac{\ln\left[(x-4)^{-1}\right]}{\ln(x)}=\frac{9}{5} \qquad (\sigma)$$

and

$$\frac{\ln\left[(2-\phi_2^4)^{-1}\right]-\ln\left[(\phi_2^4-1)^{-1}\right]}{\ln(\phi_2^4)}=2 \qquad (\tau)$$

## 2.4. NESTED RADICAL FORMS AND DIRECTED CYCLIC GRAPHS

If between three numbers $k \in \mathbf{R},\ r \in \mathbf{R},\ \alpha \in \mathbf{R},\ (k>0,\ r>1,\ \alpha>0)$ linked by the relation $\dfrac{\ln\left(\dfrac{k}{r-1}\right)}{\ln(r)}=\alpha,$ then $r = \sqrt[\alpha+1]{k+\sqrt[\alpha]{k+\sqrt[\alpha+1]{k+\sqrt[\alpha]{k+\sqrt[\alpha+1]{k+\sqrt[\alpha]{k\ldots}}}}}}$ and $k = r^{\alpha+1}-r^{\alpha}$

More generally, one real solution of the equation

$$x^\beta - \mu x^\alpha - \lambda = 0, \qquad (\beta > \alpha > 0,\ \mu > 0,\ \lambda > 0)$$

may be written as

$$\sqrt[\beta]{\lambda + \mu \sqrt[\frac{\beta}{\alpha}]{\lambda + \mu \sqrt[\frac{\beta}{\alpha}]{\lambda + \mu \sqrt[\frac{\beta}{\alpha}]{\lambda + \ldots}}}}$$

Even more generally, one real root of the equation

$$\gamma x^\beta - \mu x^\alpha - \lambda = 0, \quad (\beta > \alpha > 0, \ \gamma > 0, \ \mu > 0, \ \lambda > 0)$$

may be written as

$$\sqrt[\beta]{\frac{\lambda}{\gamma} + \frac{\mu}{\gamma} \sqrt[\frac{\beta}{\alpha}]{\frac{\lambda}{\gamma} + \frac{\mu}{\gamma} \sqrt[\frac{\beta}{\alpha}]{\frac{\lambda}{\gamma} + \frac{\mu}{\gamma} \sqrt[\frac{\beta}{\alpha}]{\frac{\lambda}{\gamma} + \frac{\mu}{\gamma} \ldots}}}}$$

It follows that the nested radical forms for $\phi_n$ with $k \in N$ and $\alpha \in N$ is

$$\phi_n = \sqrt[n+1]{1 + \sqrt[\frac{n+1}{n}]{1 + \sqrt[\frac{n+1}{n}]{1 + \sqrt[\frac{n+1}{n}]{1 \ldots}}}}$$

In particular, if $\alpha = n = 1$, then $r = \phi_1 = 1.618\ldots = \sqrt[2]{1 + \sqrt[\frac{2}{1}]{1 + \sqrt[\frac{2}{1}]{1 + \sqrt[\frac{2}{1}]{1 \ldots}}}}$

If $\alpha = n = 2$, then $r = \phi_2 = 1.465571231876\ldots = \sqrt[3]{1 + \sqrt[\frac{3}{2}]{1 + \sqrt[\frac{3}{2}]{1 + \sqrt[\frac{3}{2}]{1 \ldots}}}}$

If $\alpha = n = 3$, then $r = \phi_3 = 1.3802775\ldots = \sqrt[4]{1 + \sqrt[\frac{4}{3}]{1 + \sqrt[\frac{4}{3}]{1 + \sqrt[\frac{4}{3}]{1 \ldots}}}}$

If $\alpha = n = 4$, then $r = \phi_4 = 1.324717957244746\ldots = \sqrt[5]{1 + \sqrt[\frac{5}{4}]{1 + \sqrt[\frac{5}{4}]{1 + \sqrt[\frac{5}{4}]{1 \ldots}}}}$

(The nested *cubic root* form of the *plastic constant* – $\phi_4$ – is well known.)

### STATEMENT 3

**Self-connecting a vertex of any directed cyclic graph with *n* vertexes ($n \geq 2$) one obtains a graph whose matrix has one of its *eigenvalues* $\phi_{n-1}$**

Consider the additive sequence A005251:

1, 1, 1, 2, 4, 7, 12, 21, 37, 65, 114, 200, 351, 616, 1081, 1897, 3329, 5842… that has the recursive formula

$\forall n > 3 \quad u_n = u_{n-1} + u_{n-2} + u_{n-4}$

and the ratio limit $\phi_4^2 = 1.32471795\ldots^2 = 1.754877\ldots = \eta$ (see also page 15)

(Its *n*-th term will be designated below by $\mathbf{W}_n$.)

Now, if we write the sequence $[\eta^n]$ (the nearest integer of $\eta^n$), we obtain:

2, 3, 5, 9, 17, 29, 51, 90, 158, 277, 486, 853, 1497, 2627, 4610, 8090, 14197, … ***(a)***

which would be an additive sequence with the same recurrence equation:

$\forall n>3$  $u_n = u_{n-1} + u_{n-2} + u_{n-4}$  (just as A109377, which is, for the time being, the nearest sequence in OEIS to the sequence ***(a)***)

If we neglect the 'anomaly' in the 4$^{th}$ term (where $[\eta^4] = [9.4839092\ldots] = 9$), which 'should' have equaled 10, and if we designate by $\mathbf{V}_n$ the general term of this sequence, we immediately find that

$\forall n>4$  $V_n = P(2n+1)$  (where $P(k)$ is a Perrin number)

We also find that $\forall n>4$  $\mathbf{W}_{n+3} = \mathbf{V}_n + 2\mathbf{W}_{n-1}$

And we find too that   $\mathbf{W}_{2n} \approx \mathbf{W}_n \mathbf{V}_n$

(which, despite, the *approximate* equality, reminds an 'identity' involving Lucas and Fibonacci numbers)

## 2.5. Additive sequences, Pascal triangle, Delannoy square, generalized *p*-Tribonacci and *p*-Lucas-Tribonacci numbers

Fibonacci numbers appear in the Pascal triangle (see for, example, http://mathworld.wolfram.com/FibonacciNumber.html).

So do, in a more complicate way, the Lucas numbers (see for example http://www.mcs.surrey.ac.uk/Personal/R.Knott/Fibonacci/lucasNbs.html).

As a matter of fact, Lucas numbers appear, in a more evident way, as the sums of the shallow diagonals in the 'asymmetric additive Triangle' (that has diagonals made of 2, of odds, of squares, of sums of squares):

```
            1
          1   2
        1   3   2
      1   4   5   2
    1   5   9   7   2
  1   6  14  16   9   2
  …………………………………..
```

The diagonals of the Delannoy square

| | | | | | | |
|---|---|---|---|---|---|---|
| 1 | 1 | 1 | 1 | 1 | 1 | ……… |
| 1 | 3 | 5 | 7 | 9 | 11 | ……… |
| 1 | 5 | 13 | 25 | 41 | 61 | …….. |
| 1 | 7 | 25 | 63 | 129 | 231 | …….. |
| 1 | 9 | 41 | 129 | 321 | 681 | …….. |

………………………………………………….

sum the numbers 1, 2, 5, 12, 29, 70, etc. (Pell numbers, A000129, a sequence whose recurrence formula[12] is $\forall n>2 \ u_{n+1} = 2u_n + u_{n-1}$), while the shallow diagonals sum 1, 1, 2, 4, 7, 13, 24, 44 etc. (a sequence whose recurrence formula is $\forall n>3 \ u_{n+1} = u_n + u_{n-1} + u_{n-2}$)

Some other additive sequences appear as shallow diagonals in the Delannoy Square. There is an infinity of such diagonals and an infinity of such additive sequences:

Other shallow diagonals sum:

1, 1, 1, 2, 4, 6, 9, 15, 25, 40, 64 … (A006498, a sequence whose recurrence formula is $\forall n>3 \ u_{n+1} = u_n + u_{n-2} + u_{n-3}$)

1, 1, 1, 1, 2, 4, 6, 8, 11, 17, 27, 41, 60 … (A079972, a sequence whose recurrence formula is $\forall n>4 \ u_{n+1} = u_n + u_{n-3} + u_{n-4}$)

The Delannoy Square's shallow diagonals sum sequences with recurrence equations of the form

$\forall n>p+1 \ u_{n+1} = u_n + u_{n-p} + u_{n-(p+1)}$

and with initial conditions, 1, 1,…, 2, 4 (1 appears $p+1$ times). We'll call them (generalized) $p$-Tribonacci Numbers[13]. As $p$-Fibonacci Numbers do, they have their

---

[12] soon it will become clear that the Pell sequence is a $p$-Tribonacci sequence of order 0
[13] For $p$-Fibonacci and $p$-Lucas numbers, see also
http://www.mi.sanu.ac.yu/vismath/stakhov/index.html

(generalized) *p*-Lucas-Tribonacci (*p*-LT) companions. Their initial conditions are always:

3, 1, 1, …, 3 (the figure 1 appears exactly *p* times)

They satisfy the same recurrence relations. Besides,

$\forall n \ \forall p \quad p\text{-}LT_n = T_n + 2\,T_{n-p} + 3\,T_{n-(p+1)}$

For *p*-Tribonacci Sequences there are several combinatorial interpretations. For *p-LT* Sequences we leave the task of finding them to the reader.

## 3. RECURRENCE RELATIONS AND FORMAL GRAMMARS

**3.1.** It is possible to apply the same recursive models to linguistic operations, such as concatenation, just to take the simplest example.

Choosing as the first three words, let's say, A, AB and CA, we'll obtain

A, AB, CA, AAB, ABCA, CAAAB, AABABCA, ABCACAAAB, CAAABAABABCA, ABCACAAABCAAABAABABCA

The sequence is obviously aperiodic, but the average frequencies of each letter and perhaps of each group of consecutive letters seem to tend to a limit regardless of the concrete arbitrarily chosen words (initial conditions) and of the recursive rules. (The values of these limits of course depends on them, but their very existence not.)

This process is simple and strictly deterministic. It might be successfully (from a technical, but not necessarily esthetical[14] point of view) used in music and in composers' activity.

Let *m* be an integer at most equal to *n* and let {*i*(1), *i*(2),…, *i*(*m*)} be a part of {1, 2,…, *n*}. Assume *i*(1)< *i*(2)<…< *i*(*m*)=*n*

Choosing (in a given alphabet), as initial conditions, *n* words $W_1, W_2, …, W_n$ (distinct or not) and a recursive rule.

$\forall p > n \quad u_p = \bigoplus_{j=m}^{j=1} u_{p-i(j)}$ \hfill (A°)

(where $\oplus$ means concatenation)

---

[14] As a matter of fact, the esthetical value of a musical work depends on its author's talent, skills or genius and much less on the involved techniques. We hope to hear less cacophonic results than the words our example of sequence contains.

we indeed obtain a sequence of words with computable limits of average frequencies of every letter (in the given alphabet) and of every word of *k* letters (written in the given alphabet).

Obviously, if
  a) each of the first arbitrarily chosen words of the sequence is made of one single letter

and if

  b) these words are all distinct

then, it is possible to obtain the same sequence of words using a *context-free* formal grammar.

But if

  a') the 'initial conditions' contain words made of more than one letter

or if

  b') in the initial conditions at least one word occurs at least twice then, it is possible to obtain the same sequence of words using a formal grammar, but generally it will not anymore be *context-free*.

**3.2.** Let's slightly change the game's rules: concatenation is a non-commutative operation so we may consider some permutation of the indexes in **(10°).** It is easy to find out that recursive sequences of words constructed in this way still correspond to formal grammars that engender those words.
(For example is it easy to see that the sequence $W_1, W_2, W_3, W_2 \oplus W_1 \oplus W_3,\ldots$ whose general term will be $W_{n+3} = W_{n+1} \oplus W_n \oplus W_{n+2}$ will still generate words liable to be generated by some formal grammar.) See above the conditions under which these formal grammars are or not *context-free*.

**3.3.** Let's complicate more the rules of our game. In order to clearly explain our ideas, let's start with some preliminary definitions, conventions and notations.

Let *W* be a nonempty string, containing $n > 0$ numbers of letters. We'll write then $|W| = n$. If the set of ranks in the string at which its letters appear (namely 1, 2…, *n*) undergoes permutation, then we'll obtain what may be called a permutation of the word itself, regardless weather the letters the word is made of are all distinct or not. Choosing arbitrarily some word of length *n,* it is easy to establish a lexicographic ordering in the set of the *n*! permutations of the ranks of its letters. Suppose we have a 'computable' finite-step algorithm that permits us – for every $n$ – to chose some unique perfectly identifiable permutation.

Let's give the example of such an algorithm (further designated as **A)**: for a given $n$, we decide that the chosen permutation will always be the $m$-th permutation (in the lexicographic ordering of the set of the $n!$ permutations), setting $m$ equal to the $(n-1)^{n-1}$ –th term in the periodic sequence made of all primes smaller than $n!$. (For $n=3$, we'll have the sequence 2, 3, 5, 2, 3, 5, … etc. and the term with rank $(3-1)^{3-1}=4$ in this sequence is 2.)

As we decided that our permutation will correspond to the $(n-1)^{n-1}$-th number of this periodic sequence, then, in this specific case, we'll find the permutation 2, that is to say the second in the lexicographic order of the 3! permutations, namely the permutation (1, 3, 2). Applied to some word $W$ in the recursive sequence that we are to construct, we shall use the notation $P_{(A, |W|)}(W)$ or $P_{(A, n)}(W)$ for any word $W$ of length $|W|=n$ which underwent a permutation chosen by the algorithm **A**.

(It is easy to imagine that the sequence of words $W_1$, $W_2$, $W_3$, $W_2 \oplus P_{(A, |W(3)|)}(W_3) \oplus W_1$,….whose general term is $W_{n+3} = W_{n+1} \oplus P_{(A, |W(n+2)|)}(W_{n+2}) \oplus W_n$ can be effectively constructed.)

It seems obvious that any 'grammar' that would engender the words – and only the words – of the above described recursive sequence has infinitely many rules (even if we can express them all in a finite number of English words) and as such does not belong to 'formal' grammars. Of course, the use of these new rules at every step of the process is limited, because it depends of the steps themselves. Nevertheless, it seems that the set of the generated words is not included in the outcome vocabulary of any formal grammar. Therefore it doesn't correspond to any Turing Machine. Inasmuch it seems to be, in some sense, calculable and recursively enumerable, one can ask: how to build an automaton able to *recognize* them? How about the celebrated Church-Turing Thesis and the way it is understood or misunderstood?

### 3.4. Two open questions

1°) Given in advance the limits of the average frequencies of every letter of some given alphabet and of every word (written in the given alphabet) of at most $k$ letters (provided these limits are expressed in algebraic numbers), is it possible to obtain by an algorithm the initial conditions and the recursive rules that will generate a sequence of words with the desired average frequency limits? Is it possible always or only sometimes to find such an algorithm? If only sometimes, under which additional restrictive conditions?

2°) If the reader knows how to construct such algorithms, we would be grateful if he can also answer the question whether these algorithms are or not Polynomial-Time on the length of the input.

## 4. ABOUT A GEOMETRIC RECURSIVE PROCESS

**4.1.** One of the most beautiful recursive processes in geometry was discovered in the 19-th century by Poncelet, whose theorem is celebrated.

It was generalized by Darboux and, recently but in a different way, by Mirman, who introduced the concept of *Poncelet curve.*

Let *C* be a closed quadric and let *A* be a closed *very smooth* curve inside *C*. We conjecture that there are three and only three cases.

1°) *A* is a Poncelet-Mirman curve and then, there is at least a displacement and/or rotation of *A* which transforms it into a curve that will produce infinite never-closed polygons with an infinite set of vertexes everywhere dense in *C*.

2°) *A* gives birth to an infinite polygon that 'converges' to a finite polygon. (Its set of vertexes is attracted by a finite set of points.) Then, there is neither a displacement and/or rotation which transforms *A* into a Poncelet-Mirman curve nor a displacement and/or rotation which transforms *A* into a curve that will produce infinite never-closed polygons with an infinite set of vertexes everywhere dense in *C*.

3°) *A* is a curve that will produce infinite never-closed polygons with an infinite set of vertexes everywhere dense in *C*. Then, there is at least a displacement and/or rotation that transforms *A* into a Poncelet-Mirman curve.

## 4.2. PONCELET-MIRMAN ALGORITHMS AND GALOIS FIELDS

Boris Mirman has also found two algorithms that permit to calculate the focuses of an ellipses package generated by a 'Poncelet pair' of quadrics. Amazingly, if the number of focuses (plus the origin) is a prime, then the two matrixes that represent the order in which the focuses (and the origin) appear during computation coincide with the matrixes of a Galois field. If this number is a power of a prime, then the two matrixes coincide with the tables of a commutative ring. This fact is amazing, because these matrixes *are not* operation matrixes; they only indicate the order in which Mirman's computation method provides its outputs. However, it is hard to believe that Galois field or commutative ring structures can appear casually.

## 5. EQUIDISTRIBUTED POWER SEQUENCES

Let $R^{equi}$ designate the set of real numbers whose fractional parts of their sequence of powers is equidistributed in $]0, 1[$.

Obviously, for any close $x$ and $y$ in $R^{equi}$, frac($x^n$) and frac($y^n$) diverge apart.

**CONJECTURE**

For any ε>0 and for any $x$ and $y$ in $R^{equi}$, there is an integer $n$ such as

$|frac(x^n) - frac(y^n)| < \varepsilon$

In other words, we have chaos, except for a set of real numbers of measure 0, whose fractional parts are not equidistributed.

A similar conjecture might be formulated concerning polygons that never close in the context of would-be pairs of Poncelet-Mirman curves. Starting from two distinct arbitrarily close points the vertexes of the two polygons will diverge apart. Starting from two distinct arbitrarily chosen points on the outer quadric we'll always find on some step of the iterative process two vertexes arbitrarily close to each other. Hear again we have chaos, save for a set of cases of measure 0.

**CONCLUSION**

It seems that although recursive processes were discovered – or invented – a long time ago, they aren't yet fully investigated. It is yet not clear if generally the numbers that satisfy the equations Log[$1/(x–1)$]/Log($x$)=$k$  ($k \in N − \{1\}$) make arise not only integers but, like Pisot numbers, also *almost integers*. (Some of the numbers that satisfy the 'integer logarithmic equation' **(I)** are also Pisot numbers.)

The studied dynamical systems (and some similar) make obvious that initial pseudo-chaotic sub-sequences of sequences that turn out to converge to 0 may be as long as we wish. Their 'interpretation' is obvious. Anyway, they make arise the Feigenbaum constant. Besides, they contribute to make evident the necessity of adjusting the definition of chaos for discrete-time processes.

Recursive processes with matrixes that produces chaos are much less studied than recursive processes related to real or complex metric spaces. (Regardless of their probabilistic interpretation, Markov chain matrixes are obviously not chaotic because their row converges.)

A generalization of both recurrence (as tried in this article) and iteration (as made in the article "*Generalized iteration, Catastrophes, Generalized Sharkovsky's Ordering*" arXiv:0801.3755v2) is still to be developed by the reader interested in these issues.

ACKNOWLEDGMENTS

We express our deep gratitude to Dmitry Zotov for computer programming and to Robert Vinograd and Boris Mirman for inspiring discussions on the subject.